\documentclass[11pt]{article}

\usepackage{fullwidth}
\usepackage[top=1.5in, bottom=1.5in, left=1in, right=1in]{geometry}
\usepackage{amsgen}
\usepackage{amsmath}
\usepackage{amstext}
\usepackage{amsbsy}
\usepackage{amsopn}
\usepackage{amsfonts}
\usepackage{amssymb}
\usepackage{eepic}
\usepackage{graphicx}
\usepackage{mathrsfs}
\usepackage{faktor}
\usepackage{amsthm}

\usepackage[dvipsnames]{xcolor}
\usepackage[colorlinks, linkcolor=BrickRed,citecolor=Violet]{hyperref}
\usepackage[noabbrev,nameinlink]{cleveref}

\def\tra#1{\smash{\mathop{\mid\kern
-1pt\joinrel\relbar\joinrel\relbar}\limits^{*}_{#1}}}
\def\longtra#1{\smash{\mathop{\mid\kern
-1pt\joinrel\relbar\joinrel\relbar\joinrel\relbar}\limits^{*}_{#1}}}
\def\vlongtra#1{\smash{\mathop{\mid\kern
-1pt\joinrel\relbar\joinrel\relbar\joinrel\relbar\joinrel\relbar}\limits^{*}_{#1}}}
\def\vvlongtra#1{\smash{\mathop{\mid\kern
-1pt\joinrel\relbar\joinrel\relbar\joinrel\relbar\joinrel\relbar\joinrel\relbar}\limits^{*}_{#1}}}
\def\vvvlongtra#1{\smash{\mathop{\mid\kern
-1pt\joinrel\relbar\joinrel\relbar\joinrel\relbar\joinrel\relbar\joinrel\relbar\joinrel\relbar}\limits^{*}_{#1}}}
\def\etra#1{\smash{\mathop{\mid\kern
-1pt\joinrel\relbar\joinrel\relbar}\limits_{#1}}}

\def\bi{\begin{itemize}}
\def\ei{\end{itemize}}
\def\beq{\begin{equation}}
\def\eeq{\end{equation}}

\def\GM{\operatorname{GM}}
\def\RLM{\operatorname{RLM}}

\theoremstyle{plain}
\newtheorem{T}{Theorem}[section]
\newcommand{\bt}{\begin{T}}
\newcommand{\et}{\end{T}}
\newcommand{\ftd}{$\square$\end{T}}
\newcommand{\Eval}{\mathsf{Eval}}
\newcommand{\States}{\mathsf{States}}

\newtheorem{Proposition}[T]{Proposition}
\newcommand{\bp}{\begin{Proposition}}
\newcommand{\ep}{\end{Proposition}}
\newcommand{\fpd}{$\square$\end{Proposition}}

\newtheorem{Lemma}[T]{Lemma}
\newcommand{\bl}{\begin{Lemma}}
\newcommand{\el}{\end{Lemma}}
\newcommand{\fld}{$\square$\end{Lemma}}

\newtheorem{Corol}[T]{Corollary}
\newcommand{\bc}{\begin{Corol}}
\newcommand{\ec}{\end{Corol}}
\newcommand{\fcd}{$\square$\end{Corol}}

\newtheorem{Result}[T]{Result}
\newcommand{\br}{\begin{Result}}
\newcommand{\er}{\end{Result}}
\newcommand{\frd}{$\square$\end{Result}}

\theoremstyle{definition}
\newtheorem{Example}[T]{Example}
\newcommand{\be}{\begin{Example}}
\newcommand{\ee}{\end{Example}}

\newtheorem{Problem}[T]{Problem}
\newcommand{\bq}{\begin{Problem}}
\newcommand{\eq}{\end{Problem}}

\newtheorem{Remark}[T]{Remark}
\newcommand{\brm}{\begin{Remark}}
\newcommand{\erm}{\end{Remark}}

\newtheorem{Definition}[T]{Definition}
\newcommand{\bd}{\begin{Definition}}
\newcommand{\ed}{\end{Definition}}

\newtheorem{Construction}[T]{Construction}
\newcommand{\bco}{\begin{Construction}}
\newcommand{\eco}{\end{Construction}}

\normalbaselines
\def\abstract#1{\par\bigskip
\begingroup\small
\baselineskip=12truept
\begin{center}ABSTRACT\end{center}
\par\par\noindent
\null\hfill\hbox{\vbox{\hsize=5truein\noindent#1}}
\hfill\null\par\endgroup\par}

\usepackage[colorinlistoftodos]{todonotes}

\title{Complexity of Finite Semigroups: History and Decidability}
\author{Stuart Margolis , John Rhodes and Anne Schilling}

\date{\today}

\begin{document}

\maketitle

\begin{center}
    {\bf Abstract}
\end{center}

In recent papers Margolis, Rhodes and Schilling proved that the complexity of a finite semigroup is computable. This solved a problem that had 
been open for more than 50 years. The purpose of this paper is to survey the basic results of Krohn--Rhodes complexity of finite semigroups 
and to outline the proof of its computability.

\section{Introduction}

In two papers \cite{complexity1, complexityn}  Margolis, Rhodes and Schilling proved that the Krohn--Rhodes complexity of a finite semigroup is computable. 
This solves a problem that was open for more than 50 years~\cite{KRannals}. In this paper, we will survey Krohn--Rhodes complexity of semigroups and 
finite-state machines. From hereon in we will write ``complexity" of finite semigroups and finite-state machines. 

The purpose of this paper is to give the necessary background and an outline to read the papers~\cite{complexity1, complexityn}. We review the basic 
tools of complexity theory. We look at some lower and upper bounds to complexity. After summarizing the tools needed to prove complexity is computable 
we give an outline of the proofs in~\cite{complexity1, complexityn}.

For more background in complexity theory see~\cite[Chapters 5-9]{Arbib}, \cite[Chapters 4]{qtheory}, and~\cite{TilsonXI, TilsonXII}. See these 
references for all the definitions and results referred to in this paper. All semigroups in this paper are assumed to be finite.

A transformation semigroup (monoid) (abbreviated ts, tm) is a pair $(Q,S)$, where $Q$ is a set and $S$ is a subsemigroup (submonoid) of the 
monoid $\mathsf{PT}(Q)$ of all partial functions on $Q$ acting on the right of $Q$. We identify a semigroup $S$ with its right regular representation 
$(S^{1},S)$, where $S^{1}$ adjoins an identity if $S$ does not already have one.

We begin with the Krohn--Rhodes Theorem originally known as the Prime Decomposition Theorem for Finite Semigroups. Recall that a semigroup 
is aperiodic if all its maximal subgroups are trivial. A semigroup $S$ divides a semigroup $T$, written $S \prec T$, if $S$ is a quotient of a subsemigroup 
of $T$. There is an analogous definition of division of ts and tm. See the references above for the definition of the wreath product of ts. 

\bd
A semigroup $S$ is {\em prime} if whenever $S$ divides a wreath product $U \wr T$, then $S$ divides either $U$ or $T$.
\ed

An important semigroup that arises in the theory is the semigroup $\mathsf{RZ}(2)^{1}$ consisting of two right-zeroes and an identity. It is 
faithfully represented on a two element sets by sending the right-zeroes to the two constant functions and the identity to the identity transformation. 
This semigroup is called the flip-flop because of its connection to the theory of finite-state machines. See~\cite{Arbib} for this connection.

\bt[Krohn--Rhodes 1962]\label{KRT}
Every finite semigroup $S$ divides a wreath product of groups and aperiodic semigroups. One can choose the groups to be simple groups that 
divide $S$. We can choose the aperiodic semigroups to be the flip-flop $\mathsf{RZ}(2)^1$.

Furthermore, $S$ is prime if and only if either $S$ is a simple group or $S$ is a subsemigroup of $\mathsf{RZ}(2)^{1}$.
\et

\bc
A semigroup is aperiodic if and only if it divides a wreath product of copies of $\mathsf{RZ}(2)^{1}$.
\ec

The Krohn-Rhodes Theorem was first proved in the joint thesis of Kenneth Krohn and John Rhodes. It first appeared in print in \cite{PDT}. 
Since then it has appeared in a number of books \cite{Arbib, Eilenberg, Lallement, qtheory}. Theorem \ref{KRT} leads to the definition of 
complexity. Complexity was first defined in print in \cite{KRannals}, but had been defined and developed a number of years earlier.

\bd 
The Krohn--Rhodes complexity of a finite semigroup $S$ is the least number of groups in any decomposition of $S$ as a divisor of wreath
products of groups and aperiodic semigroups.
\ed

The Krohn--Rhodes Theorem guarantees that every finite semigroup $S$ has finite complexity. We write $Sc$ for the complexity of $S$. Clearly, we 
have $Sc=0$ if and only if $S$ is aperiodic and $Sc=1$ if $S$ is a non-trivial group. 

It is known that the full transformation semigroup on $n$-elements $T_n$ has complexity $n-1$. Therefore, if we define $C_n$ to be the set of semigroups 
of complexity at most $n$, then $C_{n}$ is properly contained in $C_{n+1}$ for all $n \geqslant 0$. We will indicate how to prove that $T_{n}c=n-1$ later 
in the paper.

The computability question can be stated as follows. Given a semigroup $S$, say by its multiplication table, the question is whether $Sc$ can be computed. 
The problem is not immediately decidable since one needs to search the infinite set of all wreath product decompositions of $S$ to find and guarantee that
 the decomposition is minimal. The main theorems of \cite{complexity1} and \cite{complexityn} prove that the computability question has a positive 
 solution.

\section{Upper and Lower Bounds to Complexity}

\subsection{Upper Bounds}

In this section, we discuss computable upper and lower bounds for complexity. Any decomposition of a finite semigroup $S$ as a divisor of a wreath 
product of groups and aperiodic semigroups gives an upper bound to $Sc$. In particular, any proof of the Krohn--Rhodes Theorem gives an upper 
bound to complexity.

The following result of Rhodes and Tilson~\cite{TilsonXI} gives an upper bound based on Green-Rees theory of finite semigroups. We assume familiarity 
with Green--Rees theory.

\bt[The Depth Decomposition Theorem]\label{Depth}
Define the depth $S\delta$ of a finite semigroup to be the longest chain of non-aperiodic $\mathcal{J}$-classes in $S$. Then $Sc \leqslant S\delta$.
\et

The proof of Theorem~\ref{Depth} gives an explicit decomposition of a semigroup $S$, where there are $S\delta$ non-trivial groups. It is well known that 
for the full transformation monoid $T_{n}$, we have $T_{n}\delta=n-1$ and thus Theorem \ref{Depth} implies that $T_{n}c \leqslant n-1$. As mentioned in 
the last section, we actually have equality in this case. On the other hand, if $\mathsf{SIS}_n$ is the symmetric inverse semigroup on $n$ elements consisting of 
all partial bijections on an $n$-set, then $\mathsf{SIS}_{n}\delta=n-1$ and  $\mathsf{SIS}_{n}c\leqslant 1$ as it is known that any inverse semigroup 
divides $S*G$, where $G$ is a group and $S$ is a semilattice. We will describe this decomposition for $\mathsf{SIS}_{n}$ in more detail in further sections.

While Theorem~\ref{Depth} depends on the length of chains of images of elements in a ts, there is another upper bound that is based on sizes of kernel 
classes of elements. Recall that if $s$ is a partial function on a set $Q$, then the kernel of $s$, denoted by $\operatorname{ker}(s)$, is the partition on 
the domain of $s$ given by $q\operatorname{ker}(s) q'$ if and only if $qs=q's$. The {\em degree} $sd$ of $s$ is the maximal cardinality of a partition 
class of $\operatorname{ker}(s)$. If $X=(Q,S)$ is a ts, then its degree $Xd$ is the maximal degree of any element $s \in S$. For example, 
$(Q,\mathsf{SIS}_{Q})d =1$ for any non-empty set $Q$. See \cite{Tilsonnumber} for a proof of the following theorem. See~\cite{cremona, deg2part2} 
for a detailed study of semigroups of degree 2.

\bt
Let $X=(Q,S)$ be a ts. Then $Sc \leqslant Xd$.
\et

This gives another proof that the symmetric inverse semigroup $\mathsf{SIS}_{Q}$ has complexity at most 1 for any set $Q$. On the other hand, if $Q$ 
is an $n$-element set, and $S$ is the semigroup of constant functions on $Q$, then $(Q,S)d=n$, but $Sc=0$.

We now look at an upper bound defined via chains of morphisms.

\bd\label{morphdef}
\mbox{}
\begin{enumerate}
\item{A morphism $f \colon S \rightarrow T$ is {\it aperiodic} if it is 1-1 on subgroups of $S$.  This is equivalent to having the inverse image of each 
idempotent in $\operatorname{Im}(S)$ being an aperiodic subsemigroup of $S$.}
\item{A morphism $S \rightarrow T$ is $\mathcal{L}'$ if, whenever $x,y$ are regular elements of $S$, then $xf=yf$ implies that 
$x\mathcal{L}y$.}
\end{enumerate}
\ed

The next theorem gives the connection between these morphisms and complexity theory. See also Remark~\ref{rlmrmk} and 
Lemma~\ref{rlmlem} below for their use in so called semi-local theory~\cite[Chapter 7-8]{Arbib} and~\cite[Chapter 4]{qtheory} as a basis for flow theory.

\bt
\mbox{}
\begin{enumerate}
\item{Let $A$ be an aperiodic semigroup and $S$ a semigroup. Then the projection $A \wr S \rightarrow S$ is an aperiodic morphism.}
\item{Let $G$ be a group and $S$ be a semigroup. Then the projection $G \wr S \rightarrow S$ is an $\mathcal{L}'$ morphism.}
\end{enumerate}
\et

We have a partial (much deeper!) converse:

\bt
\mbox{}
\begin{enumerate}
\item{{\bf The Fundamental Lemma of Complexity, Rhodes, 1969} Let $f \colon S \rightarrow T$ be a surjective aperiodic morphism. Then $Sc = Tc$.}
\item{Let $f \colon S \rightarrow T$ be a surjective $\mathcal{L}'$ morphism. Then $Sc \leqslant 1+Tc$.}
\end{enumerate}
\et

\brm
The decidability of complexity can be reduced to the question of whether $Sc=Tc$ or $Sc=1+Tc$ when $f\colon S \rightarrow T$ is a surjective 
$\mathcal{L}'$ morphism.
\erm

\bd
\mbox{}
\begin{enumerate}
\item{Let $S$ be a finite semigroup. An $Ap-\mathcal{L}'$ chain is a factorization of the trivial morphism $\tau_{S} \colon S \rightarrow 1$ of the form
$$\tau_{S}=\alpha_{0}\beta_{1}\cdots \beta_{n}\alpha_{n+1},$$
where each $\alpha_{i}$ is an aperiodic morphism and each $\beta_{j}$ is an $\mathcal{L}'$ morphism.}
\item{Let $S\theta$ be the minimal $n$ for which $S$ has an $Ap-\mathcal{L}'$ factorization.}
\end{enumerate}
\ed

\bt \label{Morchains}
Let $S$ be a finite semigroup.
\begin{enumerate}
\item{$Sc \leqslant S\theta$.}
\item{If $S$ is completely regular, then $Sc = S\theta$.  Therefore complexity is decidable for completely regular semigroups. See \cite[Chapter 9]{Arbib}.}
\item{$\mathsf{SIS}(n)\theta =n-1$.  So $\theta$ is a strict upper bound to complexity.}
\item{Let $S\widehat{\theta}= \operatorname{Min}\{T\theta \mid S$ {\text is a quotient of} $T\}.$ Then $Sc = S\widehat{\theta}.$}
\end{enumerate}
\et

We note that part 4 of Theorem \ref{Morchains} requires a search though an infinite set and thus does not prove that complexity is computable.

\subsection{Lower Bounds}

Lower bounds to complexity are more difficult to obtain. If you have a decomposition of a semigroup using $n$ groups with $n>0$, how do 
you prove that there are none that use strictly less than $n$ groups? We give one lower bound based on the concept of type I and type II 
subsemigroups and subsemigroup chains. 

\bd
Let $S$ be a finite semigroup.
\begin{enumerate}
\item{A subsemigroup $T$ of $S$ is an absolute type I subsemigroup if it is generated by a chain of $\mathcal{L}$-classes of $T$.}
\item{The type II subsemigroup $S_{II}$ of $S$ is the smallest subsemigroup of $S$ containing all idempotents and closed under weak conjugation: 
if $xyx = x$ for $x,y \in S$, then $xS_{II}y \cup yS_{II}x \subseteq S_{II}$.}
\end{enumerate}
\ed

\brm
A very important theorem of Ash~\cite{Ash} confirmed a conjecture of Rhodes and shows that if $S$ is a semigroup then $S_{II}$ consists of all 
elements of $S$ that are related to the identity element under all relational morphisms from $S$ to a group. This latter subsemigroup is also 
called the group kernel of $S$. Ash's result showing that the group kernel is precisely the type II subsemigroup of a semigroup $S$ shows that 
membership in the group kernel is decidable. This plays a crucial role in finite semigroup theory and in particular in the proof of decidability of complexity.
\erm

Here is the motivation for type I and type II subsemigroups:

\bt
Let $S$ be a semigroup.
\begin{enumerate}
\item{Assume that $S$ divides $T \wr G \wr A$ for some semigroup $T$, group $G$ and aperiodic semigroup $A$. If S' is 
an absolute type I subsemigroup of $S$, then
$S'$ divides $T \wr A' \wr G'$, where $A'$ is aperiodic and $G'$ is a group.}
\item{Assume that a semigroup $S$ divides $T \wr G$, where $T$ is a semigroup and $G$ is a group. Then $S_{II}$ divides $T$.}
\end{enumerate}
\et

\bc
Let $S$ be a semigroup with complexity $n>0$. Let $T$ be an absolute type I subsemigroup of $S$. Then $T_{II}c < Sc$.
\ec

\bt[Rhodes and Tilson 1972 \cite{lowerbounds2}] 
Let $S$ be a semigroup. Define $Sl$ to be the largest integer $k$ such that there is a chain of subsemigroups
$$ S \geqslant T_{1} > (T_{1})_{II} > T_{2} > (T_{2})_{II} > \ldots T_{k} >(T_{k})_{II},$$
where each $T_{i}$ is a non-aperiodic absolute type I semigroup and $(T_{k})_{II}$ is not aperiodic. Then $Sl \leqslant Sc$.
\et

It is known that $T_{n}$ is generated by $S_{n}$ and any element of rank $n-1$ and is thus an absolute type I semigroup. 
Furthermore $T_{n} - S_{n}$ is idempotent generated and it follows that $(T_{n})_{II} = (T_{n} - S_{n}) \cup \{1\}$. Since this semigroup 
contains a copy of $T_{n-1}$ it follows by induction that $T_{n}l= n-1 \leqslant T_{n}c$. This along with the upper bound for $T_{n}c$ that we 
obtained above from the Depth Decomposition Theorem proves that $T_{n}c=n-1$. In the case of the symmetric inverse semigroup 
$\mathsf{SIS}_{n}$, one sees that it is an absolute type I semigroup and that its type II subsemigroup is its semilattice of idempotents. 
It follows that $Sc=Sl \leqslant 1$ for all $n$.

In 1977 \cite{KernelSystems} the first example of a semigroup $S$ for which $Sl < Sc$ was published. That is, $Sl$ is in general a proper 
lower bound to $Sc$. In~\cite{TypeIIfalls}, examples of semigroups $S_{n}$ with $n>0$ such that $S_{n}l=1$ but $S_{n}c=n$ are constructed. 
These examples led to the construction of the lower bound in~\cite{HRS.2012} that turns out to be perfect for complexity 1 \cite{complexity1} and a 
modified version computes arbitrary complexity \cite{complexityn}. On the other hand, $Sl$ is the maximal local complexity function in the sense 
of \cite{localcomplex}. 

The next result shows how to compute the maximal $\mathcal{L}'$ congruence on a semigroup.

\bd
Let $S$ be a semigroup. Let $\mathscr{L}$ be the set of regular $\mathcal{L}$-classes of $S$. Let $l \in \mathscr{L},s \in S$. Then $S$ acts by 
partial functions on $\mathscr{L}$ by $l\cdot s =ls$ if $ls$ is in the same $\mathcal{D}$-class as $l$ and undefined otherwise. Let $S^{\mathcal{L}'}$ 
be the image of $S$ under this action.
\ed

\bt\label{maxLprime}
Let $S$ be a finite semigroup. Then:
\begin{enumerate}
\item{The morphism $\sigma_{L} \colon S \rightarrow S^{\mathcal{L}'}$ is a surjective $\mathcal{L}'$ morphism.}
\item{If $f \colon S \rightarrow T$ is a surjective $\mathcal{L}'$ morphism, then there is a unique surjective $\mathcal{L}'$ morphism 
$g \colon T \rightarrow S^{\mathcal{L}'}$ such that $fg=\sigma_{L}$.}
\item{The assignment of $S$ to $S^{\mathcal{L}'}$ is the object part of a functor on the category of finite semigroups.}
\end{enumerate}
\et

The type I-type II lower bound was used by Tilson and Rhodes~\cite{lowerbounds2}  to show that it is decidable if there exists an aperiodic semigroup 
$A$ and a group $G$ such that  $S$ divides  $A \wr G$. Karnofsky and Rhodes~\cite{complex1/2} used the $Ap-\mathcal{L}'$ upper bound to prove 
that it is decidable if there exists an aperiodic semigroup $A$ and a group $G$ such that  $S$ divides  $G \wr A$. Neither of these results generalize 
to the case of deciding complexity 1, that is if there exist aperiodic semigroups $A_{1},A_{2}$ and a group $G$ such that $S$ divides $A_{1}\wr G \wr A_{2}$.

\bt[Rhodes--Tilson 1972~\cite{{lowerbounds2}}]
Let $S$ be a semigroup. Then $S$ divides $A\wr G$ where $A$ is aperiodic and $G$ is a group if and only if $S_{II}$ 
is aperiodic. It follows that $Sc=Sl$. This condition is decidable.
\et

\bt[Rhodes--Karnofsky 1979~\cite{complex1/2}] Let $S$ be a semigroup. Then $S$ divides $G\wr A$, where $G$ is a group and $A$ is aperiodic if and only if 
$S^{\mathcal{L}'}$ is aperiodic. It follows that $Sc=S\theta$. This condition is decidable.
\et

\section{The Complexity Theory of Inverse Semigroups}

Inverse semigroups play an important role in their own right in many parts of mathematics. This is because they have faithful representations 
as ts of partial bijections on sets. 
In his Ph.D. thesis in 1968, Tilson proved that the symmetric inverse semigroup on a set $Q$, $\mathsf{SIS}_{Q}$ divides 
$(\{0,1\},\{0,1\}) \wr (Q,\mathsf{Sym}_{Q})$, 
the wreath product of the symmetric group on $Q$ and the two-element semilattice $\{0,1\}$. By the Preston-Wagner Theorem~\cite{Lawson},
every inverse semigroup $S$ is isomorphic to a subsemigroup of $\mathsf{SIS}_{S}$ and it follows that every inverse semigroup has complexity at most 1.

In this section, we show that many of the important concepts of inverse semigroup theory including $E$-unitary inverse semigroups and fundamental
 inverse semigroups have natural interpretations in terms of complexity theory. Congruences on ts whose quotients act by partial bijections play a crucial 
 role in complexity theory. It is useful to explicitly explain the complexity theoretic aspects of inverse semigroup theory.

\bd
We give a complexity theoretic definition of the important class of $E$-unitary semigroups. An inverse semigroup $S$ is {\em $E$-unitary} if the 
morphism $\sigma_{Gp} \colon S \rightarrow G$ from $S$ to its maximal group image is an aperiodic morphism.  This is equivalent to this morphism 
being idempotent-pure, that is, the idempotents form a congruence class of this morphism. This is one of the standard definitions of $E$-unitary semigroups.
\ed

\begin{Remark}
It follows from this definition that if $S$ is an $E$-unitary inverse semigroup, then $\tau_{S}=\sigma_{Gp}\tau_{G}$ is an $Ap-\mathcal{L}'$ 
factorization of the trivial morphism on $S$. Consequently, $S\theta \leqslant 1$. 
\end{Remark}

\bt[O'Carroll~\cite{OCarroll}]
\label{OCaroll} 
Let $S$ be an inverse semigroup. Then there is a semilattice $T$ such that $S$ is a subsemigroup of $T \wr G$, where $G$ is the maximal group 
image of $S$ if and only if $S$ is $E$-unitary.
\et

\bt[McAlister \cite{McAlisterPstuff}]
 \label{Ecov}
Let $S$ be an inverse semigroup. Then there is an $E$-unitary inverse semigroup and a surjective idempotent-separating morphism 
$T \rightarrow S$. Consequently $S\widehat{\theta}= Sc \leqslant 1$.
\et

Any inverse semigroup with 0, (for example, $\mathsf{SIS_{Q}}$) has maximal group image the trivial group. 
Thus an inverse semigroup with 0 is $E$-unitary and thus embeds into the semidirect product of a semilattice and a group by Theorem \ref{OCaroll} 
if and only if it is a semilattice. Therefore the result that every inverse semigroup divides a semidirect product of a semilattice and a group can not 
be replaced by an embedding theorem into such a semidirect product. 

Let $T$ be an $E$-unitary cover of $\mathsf{SIS}_{Q}$. Its maximal group image (remember all semigroups in this paper are finite) is isomorphic to
its minimal ideal $G$. The group $G$ maps onto the 0 of the $\mathsf{SIS}_{Q}$ by the covering morphism. We interpret the group $G$ as resolving 
the ``singularity" of the 0 of $\mathsf{SIS}_{Q}$ that is preventing the embedding of $\mathsf{SIS}_{Q}$ as a subsemigroup of a semidirect product 
of a semilattice and a group by Theorem~\ref{OCaroll}. Expanding semigroups by removing such singularities is an important tool in semigroup 
theory and in the proof of computability of complexity of finite semigroups.

We give another approach to the complexity theory of inverse semigroups via a direct product decomposition. This will generalize via the 
Presentation Lemma and the Theory of Flows to arbitrary semigroups.

Another important class of inverse semigroups is the class of fundamental inverse semigroups. An inverse semigroup is fundamental if the largest 
congruence contained in $\mathcal{H}$ is the trivial congruence. There is a maximal congruence contained in $\mathcal{H}$ on any inverse 
semigroup $S$. Its quotient is the maximal fundamental image of $S$. 

From Theorem \ref{maxLprime}  the maximal $\mathcal{L}'$ congruence on an inverse semigroup is given by the kernel of the morphism from 
$S$ to $S^{\mathcal{L}'}$. Clearly, on a regular semigroup this is the maximal congruence contained in $\mathcal{L}$. In an inverse semigroup 
we also have that this is the maximal congruence contained in $\mathcal{R}$ by inversion. Thus this is the largest congruence contained in 
$\mathcal{H}$ for an inverse semigroup. The following theorem follows immediately. This gives a complexity theoretic interpretation of 
fundamental inverse semigroups. The reader familiar with inverse semigroup can easily check that the representation defining $S^{\mathcal{L}'}$ 
is equivalent to the Munn representation~\cite{Lawson} that is the usual way to compute the fundamental image of an inverse semigroup.

\bt
Let $S$ be an inverse semigroup. Then $S^{\mathcal{L}'}$ is the maximal fundamental image of $S$. Consequently $S$ is a fundamental inverse 
semigroup if and only if $S$ is isomorphic to $S^{\mathcal{L}'}$.
\et

The next important theorem of McAlister and Reilly \cite{McAlReilly} shows that an inverse semigroup divides the direct product of a group and 
its maximal fundamental image. This is a special case of the Presentation Lemma Decomposition Theorem discussed below.

\bt[McAlister-Reilly 1975~\cite{McAlReilly}]
\label{McRe}
Let $S$ be an inverse semigroup. Then $S$ divides $(H \wr \operatorname{Sym}({\mathscr{L}})) \times S^{\mathcal{L}'}$, where $H$ is the direct 
product of the maximal subgroups of $S$, one for each $\mathcal{D}$-class, and $\mathscr{L}$ is the set of $\mathcal{L}$-classes of $S$.
\et

\section{Tools Used for the Proof of Computability of Complexity}

\subsection{Transitive Representations of Finite Semigroups}

We look at the tools and ideas used in the proof of the main theorem by reducing the problem to semigroups faithfully represented by partial 
functions acting transitively on a set.  We first review this theory.

\bd
A semigroup $S$ is right (left, bi-) transitive if $S$ has a faithful transitive right (left, right and left) action by partial functions on some set $X$.
\ed

\begin{Remark}
In the literature these are called Right Mapping, Left Mapping and Generalized Group Mapping semigroups respectively~\cite{Arbib, qtheory}.
\end{Remark}

The most important example of a transitive action is the Sch\"utzenberger representation $\operatorname{Sch}(R)$ on an $\mathcal{R}$-class 
$R$ of $S$.

For $r \in R, s \in S$, define 
\[
	\begin{cases}
      rs, & \text{if}\ rs \in R ,\\
      \text{undefined} & \text{otherwise.}
    \end{cases}
\]
There is also the dual notion of the left Sch\"utzenberger representation $\operatorname{Sch}(L)$ on an $\mathcal{L}$-class $L$ of $S$.

\begin{Lemma}
The semigroup $S$ is right (left, bi-) transitive if and only if $S$ has a unique $0$-minimal regular ideal $I(S) \approx M^{0}(A,G,B,C)$ such that 
the right (left, right and left) Sch\"utzenberger representation of $S$ is faithful on any $\mathcal{R}, (\mathcal{L}, \mathcal{R}$ and 
$\mathcal{L})$-class $R (L, R$ and $L)$ in $I(S) -\{0\}$.  
\end{Lemma}

\bd
A bi-transitive semigroup is called {\em Group Mapping} written $\operatorname{GM}$ if the group $G$ in $I(S)$ is non-trivial.
\ed

\begin{Remark}
\label{rlmrmk}
We remark the following:
\begin{enumerate}
\item{We can identify a fixed $\mathcal{R}$-class in $I(S)$ with the set $G \times B$. We obtain a faithful transformation semigroup 
$(G \times B,S)$.}
\item{The action of $S$ on $G \times B$ induces an action of $S$ on $B$.  The faithful image of this action is called the 
{\em Right Letter Mapping} (written $\RLM(S)$) image of $S$.  We thus obtain a faithful transformation semigroup $(B,\RLM(S))$.}
\end{enumerate}
\end{Remark}

\begin{Lemma}\label{rlmlem}
Let $S$ be a $\GM$ semigroup.
\begin{enumerate}
\item{$(G \times B,S)$ divides $G \wr (B,\RLM(S))$.}
\item{$S^{\mathcal{L}'} = \RLM(S)$.}
\item{For every semigroup $T$, there is a $\GM$ quotient semigroup $S$ such that $Tc=Sc$.} 
\end{enumerate}
\end{Lemma}

The following important theorem explains why flow theory is based upon $\GM$ and $\RLM$ semigroups.

\bt
The question of decidability of complexity can be reduced to the case that $S$ is a $\GM$ semigroup and whether $Sc=\RLM(S)c$ or 
$Sc = 1+\RLM(S)c$.
\et

\subsection{The Presentation Lemma and Flows}
\label{PLFlows}

The Presentation Lemma gives a necessary and sufficient condition for a $\GM$ semigroup $S$ to have the property that $Sc = \RLM(S)c$. 
The idea goes back to Tilson's proof that complexity is computable for semigroups with at most 2 non-zero $\mathcal{J}$-classes \cite{2J}. 
There are various versions in the literature. See~\cite{AHNR.1995}, \cite[Section 4.14]{qtheory} for example. 
In~\cite{HRS.2012} the notion of flows was defined, developed and proved to give an alternative way of describing the Presentation Lemma. 
We will use flows in this paper. See~\cite{FlowsI, FlowsII} for the theory of flows over aperiodic semigroups.

The best way to learn Flow Theory is to study, work through and create examples. A Master List of Examples is being maintained 
on the ArXiv~\cite{MasterList} and on the Research Gate accounts of the first and second authors. The list is dynamic and will be updated as necessary. 

Let $S$ be a $\operatorname{GM}$ semigroup with $0$-minimal ideal $M^{0}(A,G,B,C)$. The definition of flow gives a map from the state set 
of an automaton to a lattice associated with $S$. In the literature, there are two such lattices. We define these lattices and give the connections 
between them and show that they give equivalent notions of flow. The first is the set-partition lattice $\operatorname{SP}(G \times B)$. This is the lattice 
whose elements are all pairs $(Y,\Pi)$, where $Y$ is a subset of $G\times B$ and $\Pi$ is a partition on $Y$. Here $(Y,\Pi) \leqslant (Z,\Theta)$ if 
$Y \subseteq Z$ and for all $y \in Y$, the $\Pi$ class of $y$ is contained in the $\Theta$ class of $y$. 

The second lattice is the Rhodes lattice $\operatorname{Rh}_{B}(G)$. We review the basics. For more details see~\cite{AmigoDowling}. 
Let $G$ be a finite group and $B$ a finite set. A partial partition on $B$ is a partition $\Pi$ on a subset $I$ of $B$. We also consider the 
collection of all functions $F(B,G)$, $f \colon I \rightarrow G$ from subsets $I$ of $B$ to $G$. The group $G$ acts on the left of $F(B,G)$ by $(gf)(b) = gf(b)$
for $f \in F(B,G), g \in G, b \in \operatorname{Dom}(f)$. An element $\{gf \mid f \colon I \rightarrow G, I \subseteq X, g \in G\}$ of the quotient set 
$F(B,G)/G$ is called a cross-section 
with domain $I$. It should be thought of as a projectification of a cross-section of the projection from $I$ to $B$ in the usual topological sense. An
SPC (Subset, Partition, Cross-section) over $G$ is a triple $(I,\Pi,f)$, where $I$ is a subset of $B$, $\Pi$ is a partition of $I$ and $f$ is a collection of 
cross-sections one for each $\Pi$-class $\pi$ with domain $\pi$. If the classes of $\Pi$ are $\{\pi_{1},\pi_{2}, \ldots, \pi_{k}\}$, then we sometimes write 
$\{(\pi_{1},f_{1}), \ldots, (\pi_{k},f_{k})\}$, where $f_{i}$ is the cross-section associated to $\pi_{i}$. For brevity we denote this set of cross-sections by 
$[f]_{\pi}$.  We let $\operatorname{Rh}_{B}(G)$ denote the set of all SPCs on $B$ over the group $G$ union a new element 
$\Longrightarrow\Longleftarrow$  that we call {\em contradiction} and is the top element of the lattice structure on $\operatorname{Rh}_{B}(G)$.
Contradiction occurs  because the join of two SPCs need not exist. In this case we say that the contradiction is their join.

The partial order on $\operatorname{Rh}_{B}(G)$ is defined as follows. We have $(I,\pi,[f]_{\pi}) \leqslant (J,\tau,[h]_{\tau})$ if:
\begin{enumerate}
\item
{$I \subseteq J$;} 
\item
{Every block of $\pi$ is contained in a (necessarily unique) block of $\tau$;} 
\item{if the $\pi$-class $\pi_{i}$ is a subset of the $\tau$-class $\tau_{j}$, then the restriction of $h$ to $\pi_{i}$ equals $f$ restricted to $\pi_{i}$ as elements of $F(B,G)/G$. That is, $[h|_{\pi_{i}}] = [f|_{\pi_{i}}]$.}
\end{enumerate}

See \cite[Section 3]{AmigoDowling} for the definition of the lattice structure on $\operatorname{Rh}_{B}(G)$. The underlying set of the Rhodes lattice 
$\operatorname{Rh}_{B}(G)$ minus the contradiction is the set underlying the Dowling lattice on the same set and group. The Dowling lattice 
has a different partial order. For the connection between Rhodes lattices and Dowling lattices see \cite{AmigoDowling}.

We note that $\operatorname{SP}(G \times B)$ is isomorphic to the Rhodes lattice $\operatorname{Rh}_{G\times B}(1)$ of the trivial group over the 
set $G \times B$. We need only note that a cross-section to the trivial group is a partial constant function to the identity and can be omitted, leaving 
us with a set-partition pair. There are no contradictions for Rhodes lattices over the trivial group, and in this case the top element is the pair 
$(G \times B,(G \times B)^{2})$. Despite this, we prefer to use the notation $\operatorname{SP}(G \times B)$ instead of $\operatorname{Rh}_{G\times B}(1)$.

Conversely, we can find a copy of the meet-semilattice of $\operatorname{Rh}_{B}(G)$ as a meet subsemilattice of $\operatorname{SP}(G\times B)$. 
We begin with the following important definition.

\bd
A subset $X$ of $G \times B$ is a cross-section if whenever $(g,b),(h,b) \in X$, then g = h. That is, $X$ defines a cross-section of the projection 
$\theta \colon X \rightarrow B$. Equivalently, $X^{\rho} \subseteq B \times G$, the reverse of $X$, is the graph of a partial function 
$f_{X} \colon B \rightarrow G$.  An  element $(Y,\Pi) \in \operatorname{SP}(G \times B)$ is a cross-section if  every partition class $\pi$ of $\Pi$ 
is a cross-section.  
\ed

From the semigroup point of view, a cross-section is a partial transversal of the $\mathcal{H}$-classes of the distinguished $\mathcal{R}$-class, 
$R = G \times B$ of a $\GM$ semigroup. That is, the $\mathcal{H}$-classes of $R$ are indexed by $B$ and a cross-section picks at most one 
element from each $\mathcal{H}$-class. 

An element $(Y,\Pi)$ that is not a cross-section is called a {\em contradiction}. That is, $(Y,\Pi)$ is a contradiction if some $\Pi$-class $\pi$ contains 
two elements $(g,b),(h,b)$ with $g \neq h$. We note that the set of cross-sections is a meet subsemilattice of $\operatorname{SP}(G \times B)$ and 
the set of contradictions is a join subsemilattice of $\operatorname{SP}(G \times B)$.

We will identify $B$ with the subset $\{(1,b)|b \in B\}$ of $G \times B$, which as above we think of as a system of representatives of the 
$\mathcal{H}$-classes of the distinguished $\mathcal{R}$-class. Note then that $G\times B$ is the free left $G$-act on 
$B$ under the action $g(h,b)=(gh,b)$. This action extends to subsets and partitions of $G\times B$ . Thus $\operatorname{SP}(G\times B)$ 
is a left $G$-act. An element $(Y,\Pi)$ is {\em invariant} if $G(Y,\Pi) = (Y,\Pi)$. It is easy to see that $(Y,\Pi)$ is invariant if and only if:

\begin{enumerate}
 \item{$Y=G \times B'$ for some subset $B'$ of $B$.}  

 \item{For each $\Pi$-class $\pi$, $G\pi \subseteq \Pi$ and is a partition of $G\times B''$ where $B''\subseteq B'$.}
\end{enumerate}

Thus $(Y,\Pi)$ is invariant if and only if $Y=G\times B'$ for some subset $B'$ of $B$ and there is a partition $B_{1},\ldots B_{n}$ of $B'$ such that for all $\Pi$ classes $\pi$, $G\pi$ is a partition of $G\times B_{i}$ for a unique $1 \leqslant i \leq n$. 

Let $\operatorname{CS}(G\times B)$ be the set of invariant cross-sections $(Y,\Pi)$ in $\operatorname{SP}(G \times B)$. $\operatorname{CS}(G\times B)$ is a 
meet-subsemilattice of $\operatorname{SP}(G\times B)$. We claim it is isomorphic to the meet-semilattice of $\operatorname{Rh}_{B}(G)$. Indeed, 
let $(G\times B',\Pi)$ be an invariant cross-section. Using the notation above, pick  $\Pi$-classes, $\pi_{1}, \ldots , \pi_{n}$ such that $G\pi_{i}$ is a 
partition of $G\times B_{i}$ for $i=1, \ldots n$. Since $\pi_{i}$ is a cross-section its reverse is the graph of a function $f_{i} \colon B_{i}\rightarrow G$. We map 
$(G\times B',\Pi)$ to $(B',\{B_{1},\ldots , B_{n}\},[f]_{\{B_{1},\ldots , B_{n}\}}) \in \operatorname{Rh}_{B}(G)$, where the component of $f$ on $B_{i}$ is 
$f_i$. It is clear that this does not depend on the representatives $\pi$ and we have a well-defined function from 
$F \colon \operatorname{CS}(G \times B) \rightarrow \operatorname{Rh}_{B}(G)$. It is easy to see that $F$ is a morphism of meet-semilattices.

Conversely, let $(B',\Theta,[f]) \in \operatorname{Rh}_{B}(G)$. We map $(B',\Theta,[f])$ to $(G \times B',\Pi) \in \operatorname{CS}(G \times B)$,
where $\Pi$ is defined as follows. If $\theta$ is a $\Theta$-class, let $\widehat{\theta} =\{(bf,b) \mid b \in \theta\}$. Let $\Pi$ be the collection of all subsets 
of the form $g\widehat{\theta}$ where $\theta$ is a $\Theta$-class and $g \in G$. The proof that $(G \times B',\Pi) \in \operatorname{CS}(G \times B)$ 
and that this assignment is the inverse to $F$ above and gives an isomorphism between the meet-semilattice of $\operatorname{Rh}_{B}(G)$ and 
$\operatorname{CS}(G \times B)$ is straightforward and left to the reader. Furthermore, if the join of two SPC in $\operatorname{Rh}_{B}(G)$ 
is an SPC (that is it is not the contradiction) then this assignment preserves joins. Therefore we can identify the lattice $\operatorname{Rh}_{B}(G)$ 
with $\operatorname{CS}(G \times B)$ with a new element $\Longrightarrow\Longleftarrow$ added when we define the join of two elements of 
$\operatorname{CS}(G \times B)$ to be $\Longrightarrow\Longleftarrow$ if their join in $\operatorname{SP}(G\times B)$ is a contradiction. 
We use $\operatorname{CS}(G\times B)$ for this lattice as well.

We record this discussion in the following proposition.

\bp\label{Updown}
The set-partition lattice $\operatorname{SP}(G\times B)$ is isomorphic to the Rhodes lattice $\operatorname{Rh}_{G\times B}(1)$. 
The Rhodes lattice $\operatorname{Rh}_{B}(G)$ is isomorphic to the lattice $\operatorname{CS}(G\times B)$.
\ep

A {\em congruence} on a transformation semigroup $(Q,S)$ is an equivalence relation $\approx$ on $Q$ such that if $q\approx q'$ and for $s \in S$, 
and both $qs$ and $q's$ are defined, then $qs\approx q's$. Every $s \in S$ defines a partial function on $\faktor{Q}{\approx}$ by 
$[q]_{\approx}s=[q's]_{\approx}$ if  $q's$ is defined for some $q' \in [q]_{\approx}$. The quotient $\faktor{(Q,S)}{\approx}$ has states 
$\faktor{Q}{\approx}$ and semigroup $T$ which is generated by the action of all $s \in S$ on $\faktor{Q}{\approx}$. We remark that $T$ is 
not necessarily a quotient semigroup of $S$, but is in the case that $(Q,S)$ is a transformation semigroup of total functions. 

A congruence $\approx$ is called {\em injective} if every $s \in S$ defines a partial 1-1 function on $\faktor{Q}{\approx}$. It is easy to see that the 
intersection of injective congruences is injective. Therefore, there is a unique minimal injective congruence $\tau$ on any transformation semigroup 
$(Q,S)$. We call $\tau$ the Tilson congruence on a transformation semigroup because of the following proved in~\cite{Redux}. It is central to the 
theory of flows.
 
\bt\label{redux}
Let $(G \times B,S)$ be a $\GM$ transformation semigroup. Then the minimal injective congruence $\tau$ on $S$ is defined as follows: 
$(g,b) \tau (g',b')$ if and only if there are elements $s,t \in S_{II}$ such that $(g,b)s=(g',b')$ and $(g',b')t=(g,b)$.
\et
 
\brm
\mbox{}
\begin{enumerate}
\item{We can state this theorem by $(g,b)S_{II}=(g',b')S_{II}$. That is, $(g,b)$ and $(g',b')$ define the same ``right coset'' of $S_{II}$ on $G \times B$.}
\item{The proof in \cite{Redux} works on an arbitrary transitive transformation semigroup. We stated it in the case of $\GM$ transformation semigroups 
because that is how we will use it in this paper.}
\item{Theorem \ref{redux} can be used to greatly simplify the proof in \cite{lowerbounds2} for decidability of membership in $S_{II}$ for regular elements 
of an arbitrary semigroup.}
\end{enumerate}
\erm

There are three versions of the Presentation Lemma and its relation to Flow Theory in the literature \cite{AHNR.1995}, \cite[Section 4.14]{qtheory}, 
\cite{HRS.2012}. These have very different formalizations and it is not clear how to pass from one version to another. The terminology is 
different as well. For example, the definition of cross-section in each of these references, as well as the one we use in this paper is different. 
The following theorem gives a unified approach to these topics. It summarizes known results in the literature and is meant to emphasize the strong 
connections between flows, the Presentation Lemma and the corresponding direct product decomposition. For background on the derived transformation 
semigroup and the derived semigroup theorem see~\cite{Eilenberg, qtheory}. The proof of the following two theorems and their corollary
are the central results of flow theory~\cite{FlowsI}.

\bt \label{uniform}
Let $(G \times B,S)$ be $\GM$ and assume that $\RLM(S)c \leqslant n$. Then the following are equivalent:
\begin{enumerate}
\item{$Sc \leqslant n$.}
\item{There is an aperiodic relational morphism $\Theta \colon S \rightarrow H\wr T$, where $H$ is a group and $Tc \leqslant n-1$.}
\item{There is a relational morphism $\Phi \colon S \rightarrow T$, where $Tc \leqslant n-1$ and such that the Derived Transformation semigroup 
$D(\Phi)$ is in $Ap*Gp$. Equivalently, the type II subsemigroup of $D(\Phi)$ is aperiodic.}
\item{There is a relational morphism $\Phi \colon S \rightarrow T$, where $Tc \leqslant n-1$ such that the Tilson congruence $\tau$ on the 
Derived Transformation semigroup $D(\Phi)$ is a cross-section.}
\item{$(G \times B,S)$ admits a flow from a transformation semigroup $(Q,T)$ with $(Q,T)c \leqslant n-1$.}
\item{$S \prec (G \wr \operatorname{Sym}(B) \wr T) \times \RLM(S)$ for some transformation semigroup $T$ with $Tc \leqslant n-1$.}
\end{enumerate}
\et

As a special case, we restate Theorem \ref{uniform} in the case of aperiodic flows. This leads to a necessary and sufficient 
condition that a $\GM$ semigroup $S$ with $\RLM(S)c=1$ itself has complexity 1.

\bt
Let $(G \times B,S)$ be $\GM$ and assume that $\RLM(S)c \leqslant 1$. Then the following are equivalent:
\begin{enumerate}
\item{$Sc=1$.}
\item{There is an aperiodic relational morphism $\Theta \colon S \rightarrow H\wr T$, where $H$ is a group and $T$ is aperiodic.}
\item{There is a relational morphism $\Phi \colon S \rightarrow T$, where $T$ is aperiodic such that the Derived transformation semigroup 
$D(\Phi)$ is in $Ap*Gp$. Equivalently, the type II subsemigroup of $D(\Phi)$ is aperiodic.}
\item{There is a relational morphism $\Phi \colon S \rightarrow T$ where $T$ is aperiodic such that the Tilson congruence $\tau$ on the
 Derived transformation semigroup $D(\Phi)$ is a cross-section.}
\item{$(G \times B,S)$ admits an aperiodic flow.}
\item{$S \prec (G \wr \operatorname{Sym}(B) \wr T) \times \RLM(S)$ for some aperiodic semigroup $T$.}
\end{enumerate}
\et

We emphasize the direct product decomposition in Theorem~\ref{uniform}.
\bc[The Presentation Lemma Decomposition Theorem]
\label{PLDT}
Let $(G \times B,S)$ be $\GM$ and let $k > 0$.  Then $Sc \leqslant k$ if and only if $\RLM(S)c \leqslant k$ and there is a transformation semigroup 
$T$ with $Tc \leqslant k-1$ such that 
$$ (G \times B,S) \prec (G \wr \operatorname{Sym}(B) \wr T) \times (B,\RLM(S)).$$
\ec

Note that if $S$ is a $\GM$ inverse semigroup then Theorem \ref{McRe} follows from Corollary \ref{PLDT} if we take $T$ to be the trivial semigroup. 
It is immediate from the definition in fact that a $\GM$ inverse semigroup has a flow over the trivial semigroup. See \cite{FlowsI} for a description of 
semigroups that have flows over the trivial semigroup.

\subsection*{More tools used in the proof of computability of complexity 1}

The proof of computability of complexity 1 computes a finite state automaton $\mathcal{A}$ with state set $Q$ and aperiodic transformation 
semigroup $T$ such that there either there is a flow $f \colon Q \rightarrow \operatorname{SP}(G\times B)$ defined on $\mathcal{A}$ or no aperiodic flow exists.
We give a list of the tools, steps and ideas used in the proof.  

\begin{itemize}

\item{{\color{blue} Henckell's Theorem on Pointlike Aperiodic Sets} \cite{Henckell, BenHenckell}  

A subset $Y$ of a semigroup $S$ is pointlike for aperiodic semigroups if for all relational morphisms $\theta \colon S \rightarrow T$ with $T$ aperiodic, 
there is a $t \in T$ such that $Y \subseteq t\theta^{-1}$.

Henckell proved \cite{Henckell} that the pointlike subsets of $S$ are the smallest subsemigroup of the semigroup of subsets $P(S)$ containing the 
singletons, closed under subset and closed under the operation that sends a set $Z$ to $Z^{\omega}(Z \cup Z^{2} \cup \cdots \cup Z^{\omega})$,
where $Z^{\omega}$ is the unique idempotent in the subsemigroup generated by $Z$. 
\smallskip
It follows that it is decidable if a subset of $S$ is pointlike.}

\end{itemize}

\begin{itemize}
\item{{\color{blue} The Karnofsky--Rhodes Expansions}   

An {\em expansion} is a functor $F$ on the category of finite semigroups such that there is a surjective morphism $F(S) \rightarrow S$ that form part 
of a natural transformation from $F$ to the identity functor. Expansions have played a crucial role in finite semigroup theory, the first being the 
Rhodes expansion, that was used for a proof of the Fundamental Lemma of Complexity. Expansions intuitively remove ``singularities'' in semigroups 
replacing them with an expanded semigroup with better behavior. In the proof for computability of complexity, expansions are used to find a cover of 
an aperiodic semigroup whose Cayley graph look like the Cayley graph of free aperiodic Burnside semigroups. We will discuss this below.

The first expansion we look at is the Karnofsky--Rhodes expansion. We briefly recall its definition. See~\cite[Definition 4.15]{gst} 
and~\cite[Section 4.4]{complexity1} for more details and properties of the Karnofsky--Rhodes expansion. 

Let $S$ be a semigroup with a generating set $X$. The right Cayley graph $\operatorname{RCay}(S,X)$ of $S$ with respect to X is the rooted 
graph with vertex set $S^\mathcal{I}$, where $S^\mathcal{I}$ adjoins an identity $\mathcal{I}$ to $S$, even if $S$ has an identity. It has root 
$r=\mathcal{I}$ and edges $s\stackrel{a}{\rightarrow}sa$ for all $(s, a, s^{\prime}) \in S^\mathcal{I} \times X \times S^\mathcal{I}$, where 
$s^{\prime} = sa$ in $S^\mathcal{I}$. Notice that the strongly connected components of the right Cayley graph have vertices the 
$\mathcal{R}$-classes of $S^{\mathcal{I}}$. An edge $s\stackrel{a}{\rightarrow}sa$ is called a {\em transition edge} if 
$sa <_{\mathcal{R}} s$ in $S^{\mathcal{I}}$. That is, the edge goes from one connected component to one lower in the usual order on connected 
components of a directed graph. The {\em right Karnofsky--Rhodes expansion} identifies two non-empty paths starting at $\mathcal{I}$ if their 
edges multiply to the same element of $S$ and traverse the same set of transition edges. This defines a congruence on the free semigroup $X^{+}$. 
The quotient semigroup $\operatorname{RKR}(S,X)$ is called the {\em right Karnofsky--Rhodes expansion relative to $X$}. There is an obvious 
dual construction of the left Karnofsky--Rhodes expansion based on the left Cayley Graph that we denote by $\operatorname{LKR}(S,X)$.}
\end{itemize}

\begin{itemize}
\item{{\color{blue} The McCammond Expansion and the Geometric Semigroup Theory (GST) Expansion}

We will define  the McCammond expansion and apply it to the Karnofsky--Rhodes expansion. For more details on the McCammond expansion,
 see~\cite{gst} and~\cite[Section 4.4.2]{complexity1}. Recall that a simple path in $\operatorname{RKR}(S,X)$ is a path that does not visit any 
 vertex twice. Empty paths are considered simple. The McCammond expansion $\operatorname{Mc\circ RKR(S,X)}$ of $\operatorname{RKR}(S,X)$ 
 is the automaton with vertex set V, which is the set of simple paths in $\operatorname{RKR}(S,X)$. The edges are given by the following rules. If 
 $p$ is a simple path and $a \in X$ is such that $p(p,a,pa)$ is a simple path, then we have an edge from $p$ to $p(p,a,pa)$. Otherwise, there is a 
 unique vertex $v$ on $p$ such that $v=pa$. In this case, letting $q$ be the prefix of $p$ from its start vertex to $v$ we have an edge $(p,a,q)$. 

It is known that $\operatorname{Mc \circ RKR}(S,X)$ has the unique simple path property, defined as follows:  A rooted graph $(\Gamma, r)$ 
with root $r$ has the unique simple path property if for each vertex $v$ in $(\Gamma,r)$ there is a unique simple path from the root $r$ to $v$. 
The transition semigroup of $\operatorname{Mc \circ RKR}(S,X)$ is known to be aperiodic if $S$ is aperiodic~\cite{gst}.

The Geometric Semigroup Theory (GST) expansion of $(S,X)$ is the composition of the left Karnofsky--Rhodes expansion, the McCammond 
expansion and the right Karnofsky--Rhodes expansion. That is,
\[
	\operatorname{GST}(S,X) = \operatorname{RKR}\circ\operatorname{Mc}\circ \operatorname{LKR}(S,X).
\]
The $\operatorname{GST}$ expansion is used extensively in both \cite{complexity1} and \cite{complexityn}.
}
\end{itemize}

\noindent
See \cite{RS1} where both the Karnofsky--Rhodes expansion and the McCammond expansion are used in the theory of Markov Chains.

\begin{itemize}
\item{{\color{blue} Geometric properties of Cayley graphs of Free Aperiodic Semigroups}  

The Free Semigroup $F_{n}(X)$ generated by a set $X$ and index $n$ is the free semigroup in the variety defined by the identity $x^{n}=x^{n+1}$ 
and generated by $X$.

Starting with work of McCammond in 1991 \cite{Mc91}, these and other free Burnside semigroups have been extensively studied. 
For $n>2$ their word problems are decidable.

Most important for the proof of decidability of the complexity is that for $n>2$ the semigroup $F_{n}(X)$ is finite $\mathcal{J}$-above.} 
\end{itemize}

\begin{itemize}
\item{{\color{blue} Methods from Profinite Semigroup Theory} 

Profinite semigroup theory has played a central role in finite semigroups since the 1980s. See~\cite{Almeida:book} for background and applications. 
In particular, there is a proof of the type II conjecture using new properties of free profinite groups. We need further profinite results concerning 
the notion of inevitable graphs that we build into the automata and semigroups we consider in the proof of computability of complexity. 
See~\cite{complexity1, complexityn}.}
\end{itemize}

\noindent
An important application of profinite methods in complexity theory is in the theory of stable pairs \cite{Henckellstable}. See \cite{stablepairs} for a 
profinite proof of the main theorem of \cite{Henckellstable}. This result is used to prove that it is decidable if a subsemigroup of a semigroup is a 
type I subsemigroup~\cite{HRS.2012}.
    
\section{Outline of the Proof for Computability of Complexity 1}

Let $(G \times B,S)$ be $\GM$ and let $X$ be a generating set for $S$. Based on the tools we have discussed, we define a transformation 
semigroup $\Eval_{0}(S)$ whose state set $\States_0$ is a computable subset of $\operatorname{SP}(G \times B)$ and contains $(G \times B,S)$ 
as a sub-transformation semigroup. The transformation semigroup $\Eval_{0}(S)$ is a representation by partial functions of the monoid of closure 
operations on $\operatorname{SP}(G \times B)$ called the 0-flow monoid. See~\cite[Section 4]{HRS.2012}. The reader should also study the 
calculus of operations that are applied to closure operations in~\cite[Section 2]{HRS.2012}. Of particular importance is the backflow operator 
$\overleftarrow{f}$, Kleene-star $f^{*}$ and the very important loop operator $f^{\omega+*}$. Here, $f$ is a closure operator on the lattice 
$\operatorname{SP}(G \times B) \times \operatorname{SP}(G \times B)$. See~\cite[Section 5]{HRS.2012} for precise definitions. 

The main theorem of \cite{complexity1} is the following. For complexity 1, it says that the lower bound defined in~\cite[Section 5]{HRS.2012} is perfect.

\bt
\label{cequal1}
Let $(G \times B,S)$ be $\GM$. Then $Sc =1$ if and only if $\RLM(S)c \leqslant 1$ and every element of $\States_{0}$ is a cross-section.
\et

\noindent
Here are the main steps in the proof.
\begin{enumerate}
\item{Based on $(\States_{0},\Eval_{0}(S))$ we effectively define an aperiodic automaton $\mathcal{A}$ with state set $Q$. This definition involves 
all the tools we described in the previous section.}
\item{A partial flow is a partial function $F \colon Q \rightarrow \operatorname{SP}(G \times B)$ that satisfies the properties in the definition of a 
flow for all elements in $\operatorname{Dom}(F)$. We define iteratively a sequence $F_{i}$ for $i=0,1,\ldots$ of  partial flows with Domains 
$D_{i} \subseteq Q$ with $D_{i} \subseteq D_{i+1}$ that approximate a flow from $\mathcal{A}$ to the state set of $\Eval_0(S)$.}
\item{We prove that our iterative procedure either converges to a flow in a computable finite number of steps from our aperiodic automaton to 
$\operatorname{SP}(G \times B)$ if every element in $\States_0$ is a cross-section or there is no automaton with aperiodic transformation semigroup 
that defines a flow.}
\end{enumerate}

Theorem \ref{cequal1} follows immediately. Of course there are many details and the reader is invited to read~\cite{complexity1} to see all of them.

\section{Outline of the Proof for Computability of Arbitrary Complexity}

In~\cite{HRS.2012}, a computable lower bound to complexity is constructed as follows. Let $(G \times B,S)$ be a $\GM$ ts. 
A transformation semigroup $(\States_{k},\Eval_{k}(S))$ is defined for each $k \geqslant 0$. For each $k$, 
$\States_{k} \subseteq \operatorname{SP}(G \times B)$. We have $\States_{k+1} \subseteq \States_{k}$ for all $k$. The case $k=0$ was used 
in our description of the proof for complexity 1 in Theorem~\ref{cequal1}. The lower bound theorem is as follows.

\bt
\label{lowerbd}
Let $(G \times B,S)$ be a $\GM$ ts. If $\States_{k}$ contains an element that is not a cross-section, then $k+2 \leqslant Sc$. Equivalently, if $k$ is 
the least non-negative integer such that every element of $\States_{k}$ is a cross-section, then $k+1 \leqslant Sc$. 
\et

In the case $k=0$, Theorem~\ref{cequal1} says precisely that the lower bound in Theorem~\ref{lowerbd} is perfect. For $k>0$, the authors 
of~\cite{complexity1} and~\cite{complexityn}  believe that the lower bound in Theorem \ref{lowerbd} is strict.  It is conjectured 
in~\cite{complexityn} that \cite[Conjecture 1.9 and Definition 5.9]{KernelSystems} give an example of a $\GM$ semigroup $S$ with $Sc = 3$, 
but for which the lower bound of Theorem~\ref{lowerbd} is 2. This is because the notion of $k$-loopable in~\cite{HRS.2012} is probably too weak. 

In~\cite{complexityn}, we strengthen the notion of $k$-loopable in \cite{HRS.2012} precisely because of this problem.
We briefly describe our improved version of the evaluation
transformation semigroup $(\States_{k},\Eval_{k}(S))$.

Let $k \geqslant 0$. We modify the definition of the $k^{th}$-flow monoid $F_{k}(\operatorname{SP}(G \times B))$ in~\cite[Definition 4.26]{HRS.2012} 
by replacing paragraph (5) of this Definition as follows.

Let $S$ be any semigroup and let $J$ be a $\mathcal{J}$-class of $S$. Then there is a unique maximal subgroup $N_{J}^{k}$ of the maximal 
subgroup $G_{J}$ of $J$ that is pointlike with respect to the pseudovariety $C_{k}$ of semigroups of complexity at most $k$. Furthermore, 
$N_{J}^{k}$ is a normal subgroup of $G_{J}$. 

Now consider~\cite[Definition 4.26]{HRS.2012}. We replace paragraph (5) in this Definition by:
\medskip
\begin{center}
(5') $f \in M_{k}(\operatorname{SP}(G \times B))$ belongs to  $N_{J}^{k}$, then $f^{\omega+*} \in M_{k}(\operatorname{SP}(G \times B))$.
\end{center}

We still denote by $\Eval_{k}(SP(G\times B))$ the image of our modified $M_{k}(\operatorname{SP}(G \times B))$ by the representation defined 
in~\cite[Section 5]{HRS.2012} and define the set $\States_{k}$ as in that section. Then~\cite[Proposition 3.4]{complexityn} and its corollary prove by 
induction that $(\States_{k},\Eval_{k}(\operatorname{SP}(G \times B))$ is a computable ts. 

By induction we build a relational morphism from the $\operatorname{GST}$ expansion of $(\States_{k},\Eval_{k}(\operatorname{SP}(G \times B))$ 
to a ts of complexity strictly less than $k$. We then use the results of~\cite{complexity1} that the derived ts of this relational morphism has complexity 
at most 1. The Derived Semigroup Theorem~\cite{qtheory} completes the proof of the following result which is the main theorem of~\cite{complexityn}.

\bt
Let $S$ be a $\GM$ semigroup with corresponding transformation semigroup $(G \times B,S)$. Then there 
is a transformation semigroup $(Q,T)$ with 
$Tc \leqslant k - 1$ and a flow $F \colon Q \rightarrow 
\operatorname{SP}(G \times B)$  if and only if every element
of $\States_{k-1}$ is a cross-section. 
Consequently, $Sc \leqslant k$ if and only if 
$\RLM(S)c \leqslant  k$ and every element of 
$\States_{k-1}$ is a cross-section. 
\et 

Since $|\RLM(S)| < |S|$ for any $\GM$ semigroup $S$, induction on cardinality and on $k$ prove that complexity is decidable for all $k$.

We hope that this survey will help the reader to digest the proofs of computability of complexity in \cite{complexity1} and \cite{complexityn}. 
We strongly suggest that the reader read \cite{HRS.2012} along side the reading of these two papers.

\begin{center}

Acknowledgements
    
\end{center}

This paper is based on a talk by the first author at the Conference on Theoretical and Computational Algebra 2024 held in Aveiro, Portugal from 
June 30-July 5, 2024. We all thank the organizers for their kind invitation to talk and to write this paper as a summary of that lecture.

The conference also celebrated the contributions of John Meakin to mathematics. The authors thank John for his many contributions. In particular, 
the first author spent 12 years working with John at the University of Nebraska, Lincoln and they have had an active and exciting mathematical 
collaboration over many years.

The last author was partially supported by NSF grant DMS--2053350.

\bibliography{stubib.bib}
\bibliographystyle{abbrv}


\end{document}